 \documentclass[12pt]
{amsart}

\usepackage{eepic}

\advance\textheight2.8cm
\advance\voffset-1cm

\newcommand{\fbc}{^{f; \bar c}}
\newcommand{\tn}{\bar \t^{3^{n+1}}}
\def\raw{{\text{raw}}}
\def\clean{{\text{clean}}}
\def\bounded{{\text{bounded}}}
\def\constant{{\text{constant}}}

\swapnumbers

\def\xx#1 {%
\newtheorem{#1}[thm]{#1}}
\let\yy\xx

\theoremstyle{plain}
\xx Theorem
\xx Corollary
\xx Lemma
\xx {Main Lemma}
\xx {Independence Lemma}
\xx Claim
\xx Conclusion
\yy Fact
\xx {Canonization Theorem}
\xx Observation

\theoremstyle{definition}
\yy Definition
\yy Construction
\yy Notation
\yy Setup
\yy {Fact and Definition}

\theoremstyle{remark}

\yy{Case 1}
\yy{Case 2}
\yy Remark

\yy Acknowledgments
\yy Note
\yy {Overview of the construction}

\def\itm#1 {\item[(#1)]}

\newcommand{\x}{{\tt x}}
\newcommand{\y}{{\tt y}}
\newcommand{\z}{{\tt z}}
\newcommand{\s}{{\tt s}}
\renewcommand{\t}{{\tt t}}
\newcommand{\T}{{\mathbb T}}

\begin{document}

\title[No opc lattices]{There are no infinite order polynomially complete lattices after all}

\author{Martin Goldstern}
% \thanks{The first author is supported by ...}
\address{Institut f\"ur Algebra\\
Technische Universit\"at Wien\\
Wiedner Hauptstra\ss e 8-10/118.2\\
A-1040 Wien, Austria}

\email{Martin.Goldstern@tuwien.ac.at}

\author{Saharon Shelah}
\thanks{The second author is supported by the
   German-Israeli Foundation for Scientific Research \& Development
   Grant No. G-294.081.06/93.   Publication number 688.}
\address{Department of Mathematics\\
Hebrew University of Jerusalem\\
Givat Ram\\
91904 Jerusalem, Israel}
\email{shelah@math.huji.ac.il}

\subjclass{Primary 06A07; 
%       Combinatorics of partially ordered sets
%
secondary
  08A40, 
%       Operations, polynomials, primal algebras
%
  06B99,
%       None of the above but in this section (06Bxx=Lattices)
% 
03E55} 
%       Large cardinals

%%%  03E05   
%%%%    Combinatorial set theory, See also {04A20}
%%%  03E35  
%%%%    Consistency and independence results


\date{Oct 8, 1998}

\begin{abstract}
If $L$ is a  lattice with the interpolation property
 whose cardinality is a strong limit
cardinal of uncountable cofinality, then some finite power $L^n$ 
 has an antichain of size~$ \kappa $. Hence there are no infinite opc 
 lattices.   

However, the existence of strongly amorphous sets implies (in ZF) the
existence of infinite opc  lattices. 
\end{abstract}

\maketitle

% \hrule
% 
% \medskip 
% \begin{center}
% \bf {\tt preprint} ---  Please send comments to the first author
% \end{center}
% \medskip
% \hrule
% \bigskip

\setcounter{section}{-1}
\section{Introduction}

We call a lattice~$L$ {\bf $n$-order polynomially complete (opc)}  if
  every monotone  function $L^n\to L$ is 
 induced by a lattice  polynomial, 
and we say that $L$ is order polynomially complete if $L$ is $n$-order
polynomially complete for every~$n$. 

This definition is from Schweigert's Ph.D. thesis \cite{Schweigert:1972}.  
 The
survey \cite{Kaiser:1995} gives several results and bibliographical
 references for results on order polynomially complete lattices. 

\medskip
While the finite opc lattices are now well understood, the main
question on infinite opc lattices: {\em are there any?} has remained
open until now.

We showed in \cite{633}  that the size of an infinite opc\
lattice (if one exists at all) must be a strongly inaccessible
cardinal. 

We now complement this result  by showing (in ZFC) that the
cardinality of an opc lattice {\bf cannot} be a strongly
inaccessible cardinal.    Hence there are no infinite opc lattices. 

  Again the  proof 
is not algebraic in nature, but  based on a counting argument.

  Unlike our previous proof, which
employed the heavy machinery of partition calculus, this paper uses only
very basic set theory (the notions of ``cofinality'' and ``strong limit'')
and some baby model theory  (the notion of ``type''). 

We also point out that some version of AC (the axiom of choice) is
necessary for our result,  since under a strong negation of AC there are
pathological sets, which (while being infinite) are sufficiently similar
to finite sets that it is still possible to build an
 opc lattice on them.

\section{No opc lattices}

\begin{Definition} We say that a lattice $L$ has the unary
interpolation property ($1$-IP), if every monotone function from $L$
into $L$ can be interpolated by a polynomial on any finite set.
(Equivalently: Whenever $a,b,c,d\in L$,
 $b \not\le a$, $c \le d$, then there is
a polynomial $p$ with $p(a)=c$, $p(b)=d$.) 
\end{Definition}

\begin{Theorem}  Assume that $L$ has the 1-IP.  Let $\kappa $ be the
cardinality of~$L$ and assume that $\kappa$ is a strong limit cardinal
of uncountable cofinality.  Then for some natural number $n$ there is
an antichain $A \subseteq L^n$ of cardinality~$ \kappa$. 
\end{Theorem}

\def\sq#1{(#1_i: i < \kappa )}

\begin{proof}
  Let $L$ be a lattice satisfying
the assumptions of the theorem, and pick any two distinct comparable
elements of~$L$.  We will call them $0$ and~$1$, where~$0<1$.

 We will define sequences 
$\sq L$, $\sq a$, $\sq b$, $\sq {\bar c}$, $\sq n$, $\sq \tau$ such that  the
following are satisfied for all $i,j<\kappa $: 
\begin{enumerate}
\item  $\{0,1\} \subseteq L_i \subseteq L$, $|L_i|<\kappa$
\item If $i<j$, then $L_i \subseteq L_j$.
\item $a_i$ and $b_i$ are in~$L$ and realize 
 the same (quantifier-free) 
type over~$L_i$, i.e.: 
whenever $ \sigma(\x)$ and $\tau(\x)$ are unary polynomials with coefficients
in~$L_i$, then $\sigma(a_i)\le \tau(a_i)$ iff  
$\sigma(b_i)\le \tau(b_i)$.
\item $b_i \not\le a_i$
\item $n_i$ is a natural number
\item $\tau_i$ is an $n_i+1$-ary term
\item $\bar c_i = (c_i^1, \ldots, c_i^{n_i})\in L^{n_i}$
\item $\tau_i(a_i, \bar c_i) =0$, $\tau_i(b_i, \bar c_i) = 1$. 
\item For all $i<j$ we have  $a_i, b_i, c_i^1, \ldots, c_i^{n_i}\in L_j$. 
\end{enumerate}

The sequences are defined by induction on~$i<\kappa$.  In stage $i$ we
let $L_i$ be the set of ``everything'' used so far: $L_i = \{ 0,1,a_k,
b_k, c_k^\ell: k<i, \ell\le n_k\}$.   Let $ \lambda _i = \max(|L_i|,
{\aleph_0})$.   

Since every type over $L_i$ can be represented as a set of pairs of
polynomials with coefficients in~$L_i$, there are at most 
$2^{\lambda_i}$ many possible types over~$L_i$.   By our assumption,
$\lambda_i<\kappa $ implies $2^{\lambda_i} < \kappa $, so we can find two
different elements $a_i, b_i$ which have the same type over~$L_i$.
Wlog $a_i \not>b_i$. 

Now the function that maps the set $\{x: x \le a_i\}$  to $0$
and everything else (including  $b_i$) to
$1$ is monotone, so it is realized by a polynomial~$p_i(\x)$. Let $n_i$
be the number of coefficients of~$p_i$, so we can 
write $p_i(\x)$ as $\tau_i(\x,c_i^1, \ldots c_i^{n_i})$, where
$\tau_i(\x,\y_1, \ldots, \y_{n_i})$ is an $(n_i+1)$-ary term. 

This concludes the construction of our sequences.  Note that if $(k_i:
i< \kappa)$ is the increasing enumeration of a $\kappa$-size 
 subset of~$\kappa $, 
 then $(L_{k_i}:i<\kappa)$,
$(a_{k_i}: i < \kappa )$, etc., also have all the properties listed
above. 

Since $\kappa$ has uncountable cofinality, there must be some natural
number $n$ such
that  $\{i: n_i=n\}$ has cardinality~$\kappa$, so wlog (after thinning
out our sequence, if necessary) we may assume that $n_i=n$ for all
$n$.  Similarly, since there are only countably many $n+1$-ary terms,
we may assume all $\tau_i$ are equal to some fixed term~$\tau$. 

Now let $\bar d_i = (a_i, b_i, c_i^1, \ldots, c_i^n)$.  We claim that
$(\bar d_i: i< \kappa)$ is an antichain in~$L^{n+2}$.  

Indeed, pick any $i<j$ and assume that either $\bar d_i \le \bar d_j$
or 
$\bar d_i
\ge \bar d_j$.  Since $\tau$ is monotone in each argument, we either have 
$$ 0 =  \tau(a_i, \bar c_i) \le \tau(a_j, \bar c_i) \le \tau(a_j, \bar
c_j) = 0 $$
or the converse inequality, so in any case $\tau(a_j, \bar c_i)= 0$. 
Similarly we get $\tau(b_j, \bar c_i) = 1$

However, since $a_j$ and $b_j$ have the same type over~$\bar c_i, 0,
 1$, the equation $\tau(a_j, \bar c_i)=0$ implies 
$\tau(b_j, \bar c_i)=0$.   This is a contradiction, so we conclude
 that $\bar d_i$ and $\bar d_j$ are incomparable.

\end{proof}

\begin{Conclusion}\label{main}
 There is no infinite opc lattice
\end{Conclusion}
\begin{proof}
Assume that $L$ is opc.  A fortiori, $L$ has the $1$-IP. 
Let  $\kappa = |L|$.   Since $L$ is opc, we know from
 \cite{633} that $\kappa $ must be
a strongly  inaccessible cardinal, so in particular 
 $\kappa$ is a strong limit cardinal
of uncountable cofinality.  By our theorem, there is an antichain 
 $A \subseteq L$ of cardinality~$\kappa$.   But this easily  implies
that there are $2^\kappa>\kappa $ many monotone functions from $L^n$ to~$L$,
and at most $\kappa $ many  of them can be polynomials. 
\end{proof}

\section{The role of AC}

\begin{Definition}
An infinite
 set $A$ is called ``strongly amorphous'' if, for all natural numbers~$n$,
all $n$-ary relations on~$A$ are first order definable (with parameters)
 in the language of
equality.   Equivalently, $A$ is amorphous if all sets 
$R \subseteq A^n$ are in
 the Boolean algebra generated by the sets
$\{(x_1,\ldots, x_n): x_i=x_j\}$, 
$\{(x_1,\ldots, x_n): x_i=a\}$ ($a \in A$). 
\end{Definition}

While the axiom of choice (in fact, already a very weak version of AC)
clearly implies that there are no infinite strongly amorphous sets, it
is well  known  that the theory ``ZF + there is an infinite
 strongly amorphous set''
is equiconsistent with ZFC. 
   That is, ZF cannot refute the existence
of infinite strongly amorphous sets.  Hence (as we will see below), ZF
cannot refute the existence of infinite opc lattices. 

\medskip

For the rest of this section we promise not to use the axiom of choice.

\begin{Theorem}[ZF]\label{maintheorem}
 For every infinite set $L$ there is a bounded 
 lattice
$(L,\vee,\wedge,0,1,{\le})$ such that: 
\begin{quote} For all natural numbers~$n$, for all monotone functions
$f:L^n\to L$: \\ If $f$ is definable in $(L,\vee,\wedge,\le,0,1)$,
 then $f$ is induced by a polynomial. 
\end{quote}
\end{Theorem}

\begin{Remark} \label{definable}
\begin{enumerate}
\item
The ``ZF'' above means that this theorem is proved
in the usual framework of mathematics (such as given by the Zermelo
Fr\"ankel axioms for the underlying set theory) but without invoking
the axiom of choice. 
\item By ``definable'' we mean here: 
as a relation $f \subseteq L^{n+1}$, $f$ is
definable by a first order
formula (with parameters from $L$) in the language of lattice theory. 

Since all our lattices will be bounded, it will be convenient to 
include the constants (or $0$-ary operations) $0$ and $1$ into the
``language of lattice theory''.    Thus, ``definable with
parameters $c_1,\ldots, c_k$'' will mean the same as ``definable with
parameters $c_1,\ldots, c_k,0,1$''.

% $$ f(x_1,\ldots, x_n)=y \ \Leftrightarrow \exists z_1\in L\,\forall
% z_2\in L \, \cdots \,\, \varphi(c_1,\ldots, c_k, x_1,\ldots,
% x_n,y,z_1,\ldots)$$ where $\varphi$  is a Boolean combination of
% ``atomic formulas'' $\tau_1=\tau_2$, $\tau_3\le \tau_4$, \dots, and
% the $\tau_i$ are terms in the free variables $x_1,\ldots, y, z_1,
% \ldots$, possibly also using constants from~$L$.  
\end{enumerate}
\end{Remark}

\begin{Corollary}[ZF] Assume that there is an infinite strongly
amorphous set.  Then there is an infinite order polynomially complete
lattice. 
\end{Corollary}

\begin{Construction}\label{con}
Let $L$ be an infinite set, $0$ and $1$ two distinct elements of~$L$. 
 Define a lattice structure on~$L$ by requiring
 $0\le a \le 1$ for all~$a\in  L$.  
\end{Construction}

\begin{Fact} \label{fact}
Let $(L, \vee, \wedge, 0,1,\le)$ be a lattice as in
 construction \ref{con}.   Then every subset $ R \subseteq L^n$ which is
 definable
in  $(L,\vee,\wedge,0,1,\le)$
% (see  remark~\ref{definable}(2))
 with parameters $c_1, \ldots, c_k\in L$ 
 is also
 definable (with parameters  $c_1, \ldots, c_k,0,1$)  
in the language of equality [i.e, in $(L, =)$]. 
\end{Fact}

  We will abbreviate 
a situation as in fact~\ref{fact}
 by writing ``$R $ is definable from $(c_1, \ldots,
 c_k)$.''      Functions $f:L^{n} \to L$ will be treated as relations
 $f \subseteq L^{n+1} $.

\begin{Notation}
We indicate formal variables or indeterminates by a special typeface,
e.g., $\x$, $\t_1$, etc. 

We abbreviate tuples $(c_1, \ldots, c_k)$ and  $(\x_1, \ldots, \x_n)$ as
$\bar c^k$ and  $\bar \x^n$, or sometimes $\bar c$ and $\bar \x$.  We
abbreviate $(\alpha , c_1, \ldots, c_k)$ by $(\alpha, \bar c)$ or
sometimes $\alpha, \bar c$. 

 $Mon(L, L')$ is the set of all monotone maps from
 $L$ to~$L'$.   

   $\T$ is the set of all lattice-theoretic terms  in the variables $\x_1$,
 $\x_2$, \dots, $\s_1$, $\s_2$, \dots, $\t_1$, $\t_2$, \dots \
 (We include the constants $0$ and $1$ among ``lattice-theoretic terms'')
\end{Notation}

\begin{Definition} Assume that $L$, $L'$ are   isomorphic lattices, with
isomorphism $\pi: L \to L'$.  We extend $\pi $ canonically to an
isomorphism $\pi:L^n \to L^n$. For any $f:L^n \to L$ we write the
conjugate function as $\pi f$ or $\pi(f)$: 
$$ \text{for all $\bar a \in L^n$: } \ 
(\pi f)(\pi \bar a) = \pi( f(\bar a))$$
\end{Definition}

\begin{Fact}\label{constant}
If $f:L^n \to L$ is definable from $\bar c = (c_1, \ldots, c_k)$,
$\pi: L \to L$ an 
automorphism which satisfies $\pi(c_j) = c_j$ for $j=1,\ldots, k$,
then $\pi f = f$. 
\end{Fact}

% 
% 
%  Let $g:L \to \T$ be definable from $\bar c=\{c_1,
% \ldots, c_k\}$, $A = L \setminus \{ c_1, \ldots, c_k\}$.   Then $g$ is
% constant on $A$. 
% \end{Fact}
% 
% \begin{proof} For any $a,a'\in A$ there is an automorphism $\pi$ of $L$
% which fixes $\{c_1,\ldots, c_k\}$ pointwise.   The definition of $f$
% is invariant under this automorphism, so $f = f \circ \pi$. 
% Hence $f(a) = f(\pi(a))=f(a')$. 
% \end{proof}

  \begin{Definition} Let $(L, \le)$ be a partial order, $A \subseteq
  L$. The ``monotone characteristic function of~$A$'' is the function
  $\chi_A:L \to \{0,1\}$ defined by 
  $$ \chi_A(x) =  \left\{ 
  	\begin{array}{ll}
  		1  & \mbox{if $a\ge b$ for some $b \in A$}\\
  		0  & \mbox{otherwise,}
  	\end{array}
  \right.
  $$
  i.e., $\chi_A$ is the customary characteristic function of the upward
  closure of~$A$. 
  \end{Definition}
  
  {}From now on $L$ will be a lattice as in \ref{con}
 
 \begin{Definition} 
 Let $\{c_1, \ldots, c_k\} \subseteq L$, $\bar d = (d_1, \ldots, d_m)
 \in L^m$.  We say that $\bar d$ is ``independent over~$\{c_1,\ldots,
 c_k\}$'' (or: over~$(c_1, \ldots, c_k)$) iff 
 \begin{quote} all $d_i$ are distinct, and no $d_i$ is in~$\{c_1,
 \ldots, c_k,0,1\}$. 
 \end{quote}
[This is a special case of the usual model-theoretic notion of
 independence.] 
 \end{Definition}

% Since we use a well-ordered set of variables, the set $\T$ is
% canonically well-ordered. 

  \begin{Definition}\label{chidef}
\begin{enumerate}
\item Let $\chi(\s, \x, \t_1, \t_2, \t_3) $ be the term 
$$ 
[ (\x \wedge \s) \vee \t_1] \wedge [ (\x \wedge \s) \vee \t_2]. $$ 
\item For $\alpha \in L \setminus \{0,1\}$ let
 $\chi_\alpha (\x, \t_1,\t_2, \t_3) = \chi(\alpha, \x,  \t_1,\t_2,
 \t_3)$. 
\item If $A \subseteq \{c_1 \ldots, c_k,0,1\} \subseteq L$, then 
  we define $\chi_A^{\bar c}(\bar \s^k,\x, 	\bar\t^3)$  
  (the
  ``monotone 
  characteristic function of~$A$, given the parameters $\bar c$'') as
  follows:
  \begin{enumerate} 
  \item If $A = \emptyset$, then
   $\chi_A^{\bar c}(\bar \s,\x, \t_1, \t_2,\t_3) $ is the constant term~$0$. 
  \item If $A = \{1\}$, then $\chi_A^{\bar c}(\bar \s,\x, \t_1, \t_2,\t_3) = 
  \mu( \x \wedge \t_1, \x \wedge \t_2, \x \wedge \t_3)$, where $\mu$ is
 the following ``majority term'': 
 $$ \mu(\x,\y,\z) = (\y \vee \z) \wedge (\z \vee \x) \wedge (\x \vee
 \y)$$ 
 \item If $0 \in A$, then 
  $\chi_A^{\bar c}(\bar \s,\x, \t_1, \t_2,\t_3) $ is the constant term
~$1$. 
 \item Otherwise we let $I = I_A^{\bar c} = \{i: c_i \in A \setminus
 \{0,1\}\}$ and we let $$\chi_A^{\bar c}(\x, \bar \t) = \bigvee_{i \in I}
 \chi_{c_i}(\x, \bar \t)$$
  \end{enumerate}
\item 
 If $A$ is cofinite, 
  $L \setminus A \subseteq \{c_1, \ldots, c_k\}$, then we define
 	$\chi_A^{\bar c}(\bar \s,\x, 	\bar\t)$  similarly (dually),
         such that 	fact \ref{chi} below holds.   We leave the
  details to the  	reader. 
\end{enumerate}
 \end{Definition}
 
 \begin{Fact}\label{chi} 
 If $A \subseteq \{c_1, \ldots, c_k,0,1\}$ or $L \setminus A \subseteq
  \{c_1, \ldots, c_k,0,1\}$, 
  and $\bar d = (d_1,d_2,d_3)$ is  independent over~$c_1,\ldots,
  c_k$, then the function 
 $$ a \mapsto  \chi_A^{\bar c}(\bar c, a, d_1, d_2, d_3)$$
 is the monotone characteristic function of~$A$. 
 \end{Fact}
 \begin{proof} Easy computation.  \end{proof}

 \begin{Lemma}\label{lemma}
  Let $L$ be as in construction \ref{con}. 
 Then there is a function 
 $$p: 
 \bigl( \bigcup_{n=0}^\infty Mon(L^n,L)\bigr)
\times 
\bigl(\bigcup_{k=1}^\infty L^k \bigr) 
 \ \to \T$$ that
 assigns to each pair $(f,\bar c)$ a polynomial $p\fbc\in \T$
  such
 that the following hold: 
 \begin{enumerate}
 \item If $\bar c = (c_1, \ldots, c_k)$ and $f: L^n \to L$, then 
 $p\fbc $ is a term using (at most)  the variables 
 $\s_1, \ldots, \s_k, \x_1, \ldots, \x_n, \t_1, \ldots, \t_{3^n}$.
 We will write 
 $p\fbc $ as $p\fbc(\bar \s; \bar \x; \bar \t)$. 
 \item If $\bar c$ and $ f$ are as above, and $f$ is  definable
from  $ c_1, \ldots, c_k$, then we have: 
 \begin{quote} Whenever $\bar d = (d_1, \ldots d_{3^n})$ is independent
 over  $\{c_1, \ldots, c_k\}$, \\then
 for all $\bar a \in L^n$: $ f(\bar a) = p\fbc(\bar c; \bar a;
 \bar d)$, 
 \end{quote}
that is, the function $f$ is induced by the lattice polynomial
 $p \fbc(\bar c; \bar\x; \bar d)$. 
\item Moreover, the term $p\fbc$ depends only on the ``isomorphism
 type'' of $(f, \bar c)$.   That is, whenever $\pi: L \to L'$ is an
 isomorphism, then $p\fbc(\bar \s,\bar\x,\bar \t)
 = p ^{\pi f, \pi(\bar c)} (\bar \s,\bar\x,\bar \t) $. 
 \end{enumerate}
 \end{Lemma}

\begin{Remark} Let $L$ be as above.  Then  the following 
is {\em not} provable in ZF: 
\begin{enumerate}
\item [$(*)$] There is a map $p$ which assigns to each definable monotone
function $f:L \to L$ a polynomial $p^f(\x)$ with coefficients in $L$
such that $p^f( a)= f( a)$ for all $a\in L$.
% \item $L$ is Dedekind-infinite, i.e., there is a 1-1 map from the
%        natural numbers into $L$. 
\end{enumerate}
This explains why we have to explicitly mention the parameters $\bar c$
in lemma~\ref{lemma}. 
\end{Remark}

 \begin{proof}[Proof of lemma~\ref{lemma}]
  We define the map $p$ by induction on~$n$ (the arity of
 $f$).  We leave the case $n=0$ to the reader.

  Let $f:L^{n+1} \to L$, $\bar c = (c_1,\ldots, c_k)\in L^k$. 
(Wlog assume that $f $ is definable from $\bar c$, otherwise set $p\fbc=0$.) 

Let $A = L \setminus \{c_1,\ldots,c_k, 0,1\}$. 

We will now apply the induction hypothesis to the functions
$f_\alpha:L^n\to L$, where
$f_\alpha(\bar a) = f(\bar a,\alpha)$.
  For each $\alpha\in L$ we thus get
a term $\tau_\alpha = p^{f_\alpha}$. 
 We will show that only finitely many different terms
actually appear, and that they can be combined to yield a term $p^f$.
The cases $\alpha\in A$ and $\alpha\notin A$ have to be treated
in different ways.

To save us some cases distinctions,  we agree that $c_{k+1}=0$, $c_{k+2}=1$.

For each $\alpha \in \{c_1,\ldots, c_{k+2}\}$ 
the function $f_\alpha :L^n\to L$, $f_\alpha (\bar a)
= f(\bar a, \alpha )$ is definable from $\bar c=(c_1,\ldots, c_k)$. Let 
$$ \tau\fbc_\alpha (\bar \x^n; \bar \s^k; \tn) = 
p^{f_\alpha , \bar c} (\bar \x^n; \bar \s^k; \bar \t^{3^n}) $$
(so the indeterminates  $\t_{3^n+1}, \ldots, \t_{3^{n+1}}$ do not appear in
$\tau_\alpha^{f,\bar c}$.)
By our inductive assumption we know for all $\bar a \in L^n$: 
\begin{equation}\label{1}
 \tau\fbc_\alpha (\bar a^n; \bar c^k; \bar d^{3^{n+1}}) =
f_\alpha (\bar a) = 
 f(\bar a,\alpha )
\end{equation}
for all  $\alpha \in
\{c_1,\ldots, c_{m+2}\}$, whenever $\bar d$ is independent over~$\bar c$.

Let $$\tau\fbc_\constant(\bar\x^n,\x_{n+1};\bar s^k;\bar d^{3^{n+1}})
= \bigvee_{i=1}^{m+2}\biggl[
          \chi_{c_i}(\x_{n+1},\bar \t^3 ) 
           \wedge 
          \tau\fbc_{c_i}(\bar \x^n;\bar\s^k;\tn)\biggr] .  $$ 
Hence  we have for all $\bar a \in L^n$:
\begin{equation}\label{2}
\tau_\constant(\bar a^n, \alpha; \bar c^k;\bar d^{3^{n+1}}) = 
    \begin{cases} f(\bar a, \alpha) & \text { if }
			 \alpha \in \{c_1, \ldots, c_k, 0,1\}\\
		  f(\bar a, 0)     & \text { if }
			\alpha \in A = L \setminus \{c_1,\ldots, c_k,0,1\} .
    \end{cases}
\end{equation}
whenever $\bar d $ is independent over $\bar c$.
  This follows easily from equation \eqref{1} and the monotonicity of
~$f$. 

% 
% [Proof of equation~\eqref{2}: 		
% 
% \begin{align*}
% \tau_\constant(\bar a^n, \alpha; \bar c^k;\bar d^{3^{n+1}}) &= 
%  \bigvee_{i=1}^{m+2}\bigl(
%           \chi_{c_i}(\alpha ,\bar d^3 ) 
%            \wedge 
%           \tau\fbc_{c_i}(\bar a^n;\bar c^k;\bar d)\bigr)  =\\
% &=  \bigvee_{i=1}^{m+2}\bigl(
%           \chi_{c_i}(\alpha ,\bar d^3 ) 
%            \wedge 
% 	f(\bar a, c_i) \bigr) 
% \end{align*}
% 
% If $\alpha = 1$, then this expression yields 
% $  \bigvee_{i=1}^{m+2} 	f(\bar a, c_i)$, and 
%  if $\alpha = c_j \not=1$, we get $f(\bar a, 0) \vee f(\bar a,
% c_j)$,
% and if $\alpha \in A$, we get $f(\bar a, 0)$.  The monotonicity of~$f$
% now implies equation~\eqref{2}.] 

\bigskip

For each $\alpha \in A$ the function $f_\alpha: L^n\to L$,
$f_\alpha(\bar a^n) = f(\bar a^n, \alpha)$ is definable from $(\alpha,
\bar c)$. Hence the polynomial 
$ p^{f_\alpha;\alpha,\bar c}(\bar \x^n; \x_{n+1},\bar \s^k;\bar
t^{3^n})$
satisfies for all $\bar a \in L^n$
\begin{equation}\label{3}
 p^{f_\alpha;\alpha,\bar c}(\bar a; \alpha ,\bar c^k;\bar d) =
f(\bar a, \alpha)
\end{equation}
for all $ \alpha \in A$, 
whenever $\bar d$ is independent over  $\alpha , \bar c$. 

We claim that  the function $ \alpha \mapsto 
 p^{f_\alpha;\alpha,\bar c}(\bar \x; \x_{n+1},\bar \s;\bar \t)$ is
 actually  constant on~$A$.  The reason is that $A$ is disjoint to the
 set of parameters from which $f$ is defined.   

{ \small 
  More formally, let $\alpha, \beta \in A$.  We will show 
$p^{f_\alpha;\alpha,\bar c} = p^{f_\beta ; \beta,\bar c}$. Let $\pi:L \to
L$ be an automorphism that fixes $\{c_1,\ldots, c_k,0,1\}$ pointwise, and
maps $\alpha $ to~$\beta$. Then $\pi f = f$, by
fact~\ref{constant}.

 Hence
\begin{align*}
 p^{f_\alpha; \alpha, \bar c} &= p^{\pi(f_\alpha);\pi(\alpha),\pi(\bar c)} 
&& \text{by induction hypothesis, lemma~\ref{lemma}(3)}\\
& = p^{(\pi f)_{\pi( \alpha)}; \beta, \bar c } 
&& \text{since  $\pi(f_\alpha) = (\pi f)_{\pi(\alpha)}$}
\\
& = p^{ f_\beta; \beta, \bar c } 
&&\text{since $\pi f = f $}
\end{align*}

}

    We will write 
$\tau_\raw\fbc(\bar \x, \x_{n+1};\bar \s;\bar \t)$ for the  common
 value of this function.  So we can rewrite equation \eqref{3} as: 
for all $\bar a\in L^n$, $\alpha \in A$:
\begin{equation}\label{3.5}
 \tau\fbc_\raw(\bar a, \alpha; \bar c; \bar d) = f(\bar a , \alpha)
\end{equation}
whenever $\bar d$ is independent over $\bar c$. 

The restriction that $\bar d$ has to be independent not only of $\bar c$
but also of $\alpha$ is inconvenient.   We get rid of it with the
following 
``error correction'' device: Let 
 \begin{itemize}
 \item [] $\sigma\fbc_1(\bar \x^n,\x_{n+1};\bar \s;\t_1,\ldots,\t_{3^{n+1} }) = 
 \tau_\raw\fbc(\bar \x^n,\x_{n+1};\bar \s;\bar \t_1,\ldots, \t_{3^n})$
 \item [] $\sigma\fbc_2(\bar \x^n,\x_{n+1};\bar \s;\t_1,\ldots,\t_{3^{n+1} }) = 
 \tau_\raw\fbc(\bar \x^n,\x_{n+1};\bar \s;\t_{3^n+1}, \ldots, \t_{2\cdot3^n})$
 \item [] $\sigma\fbc_3(\bar \x^n,\x_{n+1};\bar \s;\t_1,\ldots,\t_{3^{n+1} }) = 
 \tau_\raw\fbc(\bar \x^n,\x_{n+1};\bar \s; \t_{2\cdot 3^n+1}, \ldots,
\t_{3\cdot3^n})$ 
 \end{itemize}
and let $$\tau\fbc_\clean
 = \mu(\sigma_1\fbc, \sigma_2\fbc, \sigma_3\fbc),$$
 where
$\mu(\x,\y,\z) $ is a majority term,  i.e., $\mu$ satisfies 
\begin{equation}\label{majority}
\mu(a,a,b)=\mu(a,b,a)=\mu(b,a,a)=a
\end{equation}
for all $a,b\in L$, see definition~\ref{chidef}. 

Let $$\tau\fbc_\bounded(\bar\x^n, \x;\bar \s; \bar \t) =
\tau\fbc_\clean(\bar\x^n, \x;\bar \s; \bar \t)  \wedge 
 \tau\fbc_1(\bar\x^n;\bar \s; \bar \t) \wedge \chi^{\bar c}_A(\x,\bar \t^3).$$

We claim  that 
$\tau\fbc_\bounded(\bar \x^n,\x_{n+1}; \bar \s^k;\tn)$
satisfies 
\begin{equation}\label{4}
 \tau\fbc_\bounded(\bar a^{n},\alpha ; \bar c^k;\bar d^{3^{n+1}}) =
\begin{cases}
f(\bar a^n, \alpha) & \alpha \in A
\\
0 & \alpha \in \{c_1, \ldots, c_k, 0\}\\
\le f(\bar a, 1)  &  \alpha =1 
\end{cases}
\end{equation}
for any  $\bar a \in L^n$, $\alpha \in L$, 
 whenever $\bar d$ is independent over $\bar c$. 

This is clear if $ \alpha \in \{c_1,\ldots, c_k, 0\}$.    For $\alpha
= 1$ recall that $\tau\fbc_1(\bar a;\bar c;\bar d) = f(\bar a, 1)$. 

Finally,  let $\alpha \in A$, and assume that  $\bar d$  is
independent over $\bar c$. Let 
 $$
  \bar d_1 = (d_1, \ldots, d_{3^n}) \ \ \ 
 \bar d_2 = (d_{3^n+1}, \ldots, d_{2\cdot 3^n}) \ \ \ 
 \bar d_3 = (d_{2\cdot3^n+1}, \ldots, d_{3\cdot 3^n})
 $$
 Now since all the $d_j$ are distinct, 
 at least two among $\bar d_1$, $\bar d_2$, $\bar d_3$ do not contain
 $\alpha $ and hence are independent over~$\{c_1, \ldots, c_k, \alpha\}$.
 So by equation~\eqref{3.5}, at least two of the equations 
 $$ f_\alpha(\bar a) = {\sigma} \fbc_\ell(\bar a, \alpha; \bar c;\bar
 d^{3^{n+1}}) 
 \ \ \ \ \ \ell=1,2,3$$
 are true, hence we have (by the  property \eqref{majority}) 
\begin{equation*}
 \tau\fbc_\clean(\bar a, \alpha; \bar c;\bar d)
= 
f_\alpha(\bar a) 
\end{equation*}
so 
$ \tau\fbc_\bounded (\bar a, \alpha; \bar c;\bar d)
= 
f_\alpha(\bar a) \wedge \tau_1(\bar a; \bar c;\bar d) = 
f(\bar a, \alpha ) \wedge f(\bar a, 1) = f(\bar a, \alpha )$ for all
$\alpha \in A$. 

This concludes our discussion of the term $\tau\fbc_\bounded$.

\smallskip
We can now  define $p\fbc$ as 
$$ p\fbc(\bar \x^n, \x_{n+1} ; \bar \s^k; \tn) = 
\tau\fbc_\constant(\bar \x^n,\x_{n+1} ; \bar \s; \bar\t) 
\ \vee \ 
\tau\fbc_\bounded(\bar \x^n,\x_{n+1}  ; \bar \s; \bar\t) 
$$
{}From \eqref{2} and \eqref{4} we now get the desired property 
$$ p\fbc(\bar a,\alpha  ; \bar c;\bar d) =  f(\bar a, \alpha)$$
whenever $\bar d $ is independent over $\bar c$. 
 \end{proof}

Theorem \ref{maintheorem} now immediately follows from
lemma~\ref{lemma}.

% 
% To check that this works, let $\bar d^{3^{n+1}}$ be independent over
% $\bar c$. 
% 
% For $\alpha \in \{c_1, \ldots, c_k, 0, 1\}$ we have 
% $p\fbc(\bar a, \alpha; \bar c; \bar d) = f(\bar a, \alpha) \vee 0 =
% f(\bar a, \alpha)$, and for 
% $\alpha \in A$ we have 
% $p\fbc(\bar a, \alpha; \bar c; \bar d) = f(\bar a,0 ) = 
% f(\bar a, \alpha)$.

\bibliographystyle{plain}
\bibliography{%listb,
other,goldstrn}
\end{document}